\documentclass[12pt]{article}
\usepackage{amsmath,amsfonts,amssymb,amsthm}
\input amssym.def
\begin{document}
\def \Z{\Bbb Z}
\def \C{\Bbb C}
\def \R{\Bbb R}
\def \Q{\Bbb Q}
\def \N{\Bbb N}

\def \A{{\mathcal{A}}}
\def \D{{\mathcal{D}}}
\def \E{{\mathcal{E}}}
\def \H{\mathcal{H}}
\def \S{{\mathcal{S}}}
\def \wt{{\rm wt}}
\def \tr{{\rm tr}}
\def \span{{\rm span}}
\def \Res{{\rm Res}}
\def \Der{{\rm Der}}
\def \End{{\rm End}}
\def \Ind {{\rm Ind}}
\def \Irr {{\rm Irr}}
\def \Aut{{\rm Aut}}
\def \GL{{\rm GL}}
\def \Hom{{\rm Hom}}
\def \mod{{\rm mod}}
\def \ann{{\rm Ann}}
\def \ad{{\rm ad}}
\def \rank{{\rm rank}\;}
\def \<{\langle}
\def \>{\rangle}

\def \g{{\frak{g}}}
\def \h{{\hbar}}
\def \k{{\frak{k}}}
\def \sl{{\frak{sl}}}
\def \gl{{\frak{gl}}}

\def \be{\begin{equation}\label}
\def \ee{\end{equation}}
\def \bex{\begin{example}\label}
\def \eex{\end{example}}
\def \bl{\begin{lem}\label}
\def \el{\end{lem}}
\def \bt{\begin{thm}\label}
\def \et{\end{thm}}
\def \bp{\begin{prop}\label}
\def \ep{\end{prop}}
\def \br{\begin{rem}\label}
\def \er{\end{rem}}
\def \bc{\begin{coro}\label}
\def \ec{\end{coro}}
\def \bd{\begin{de}\label}
\def \ed{\end{de}}

\newcommand{\m}{\bf m}
\newcommand{\n}{\bf n}
\newcommand{\nno}{\nonumber}
\newcommand{\nord}{\mbox{\scriptsize ${\circ\atop\circ}$}}
\newtheorem{thm}{Theorem}[section]
\newtheorem{prop}[thm]{Proposition}
\newtheorem{coro}[thm]{Corollary}
\newtheorem{conj}[thm]{Conjecture}
\newtheorem{example}[thm]{Example}
\newtheorem{lem}[thm]{Lemma}
\newtheorem{rem}[thm]{Remark}
\newtheorem{de}[thm]{Definition}
\newtheorem{hy}[thm]{Hypothesis}
\makeatletter \@addtoreset{equation}{section}
\def\theequation{\thesection.\arabic{equation}}
\makeatother \makeatletter

\begin{center}
{\Large \bf  Associating quantum vertex algebras to Lie algebra $\gl_{\infty}$ }
\end{center}

\begin{center}
{Cuipo Jiang\footnote{Supported by China NSF grant 10931006,  China RFDP grant 2010007310052, and the Innovation Program of Shanghai Municipal Education Commission (11ZZ18)}\\
Department of Mathematics\\
Shanghai Jiaotong University, Shanghai 200240, China\\
Haisheng Li\footnote{Partially supported by NSA grant
H98230-11-1-0161 and China NSF grant (No. 11128103)}\\
Department of Mathematical Sciences\\
Rutgers University, Camden, NJ 08102, USA}
\end{center}

\begin{abstract}
In this paper, we present a canonical association of quantum vertex
algebras and their $\phi$-coordinated modules to Lie algebra
$\gl_{\infty}$ and its $1$-dimensional central extension. To this
end we construct and make use of another closely related infinite-dimensional Lie
algebra.
\end{abstract}

\section{Introduction}

It has been known that vertex algebras can be canonically associated to both
twisted and untwisted affine Lie algebras. More specifically, for an untwisted affine Lie
algebra $\hat{\g}$ and for any complex number $\ell$, one has a
vertex algebra $V_{\hat{\g}}(\ell,0)$ (cf. \cite{fz}, \cite{dl}, \cite{lian}, \cite{li-local}), based on a distinguished  level $\ell$
generalized Verma $\hat{\g}$-module
often called the vacuum module, while the category of
$V_{\hat{\g}}(\ell,0)$-modules is canonically isomorphic to the
category of restricted $\hat{\g}$-modules of level $\ell$. On the
other hand, it was known (see \cite{li-twisted}; cf. \cite{flm}) that the category of restricted modules for a twisted affine
algebra $\hat{\g}[\sigma]$ with $\sigma$ a finite-order automorphism
of $\g$ is canonically isomorphic to the category of
$\bar{\sigma}$-twisted modules for the vertex algebra
$V_{\hat{\g}}(\ell,0)$ (which was associated to the {\em untwisted} affine
algebra), where $\bar{\sigma}$ is the corresponding automorphism of
vertex algebra $V_{\hat{\g}}(\ell,0)$.

In this paper, we study Lie algebra $\gl_{\infty}$ in the content of quantum vertex
algebras in the sense of \cite{li-qva1}, and
as the main result we obtain a canonical association of quantum vertex
algebras to the one-dimensional central extension $\widetilde{\gl}_{\infty}$
of Lie algebra $\gl_{\infty}$, which is somewhat similar to
the association of vertex algebras and twisted modules to twisted affine
Lie algebras.

By definition, $\gl_{\infty}$ is
the Lie algebra of doubly infinite complex matrices with only
finitely many nonzero entries under the usual commutator bracket.
 A canonical
base consists of matrices $E_{m,n}$ for $m,n\in \Z$, where $E_{m,n}$
denotes the matrix whose only nonzero entry is the $(m,n)$-entry
which is $1$. For the natural representation of $\gl_{\infty}$ on
$\C^{\infty}$ with the standard base denoted by $\{v_{n}\ |\ n\in
\Z\}$, we have
$$E_{m,n}v_{r}=\delta_{n,r}v_{m}\ \ \ \mbox{  for }m,n,r\in \Z.$$
On Lie algebra $\gl_{\infty}$, there is a $2$-cocycle $\psi$ defined by
\begin{eqnarray*}
&&\psi(E_{i,j},E_{j,i})=1=-\psi(E_{j,i},E_{i,j})
\ \ \ \mbox{if $i\le 0$ and $j\ge 1$},\nonumber\\
 &&\psi(E_{i,j},E_{m,n})=0 \ \ \ \mbox{otherwise}.
\end{eqnarray*}
Using this $2$-cocycle one obtains a $1$-dimensional central
extension
$$\widetilde{\gl}_{\infty}=\gl_{\infty}\oplus \C {\bf k}.$$

With various motivations, one of us (H.L) has extensively studied
(see [Li3-5]) vertex algebra-like structures generated by fields (on
a general vector space), that behave well, but are not necessarily
mutually local (cf. \cite{dl}, \cite{li-local}). (A result of
\cite{li-local} is that every set of mutually local fields on a
general vector space canonically generates a vertex algebra.) In
this study, a theory of (weak) quantum vertex algebras and their
modules was developed, where the notion of quantum vertex algebra
naturally generalizes the notion of vertex algebra and that of
vertex super-algebra. In this theory, a key role is played by the
notion of $\S$-locality due to Etingof-Kazhdan, which is a
generalization of that of locality, and the essence is that every
$\S$-local subset of fields on a vector space $W$ generates a (weak)
quantum vertex algebra with $W$ as a natural module. (The pioneer
work \cite{ek} has been an important inspiration for the development
of this theory.)

We next explain how quantum vertex algebras are associated to Lie algebra $\widetilde{\gl}_{\infty}$.
Note that in the association of vertex algebras to affine Lie algebras,
a key role was played by the canonical generating functions (or fields), where
roughly speaking, the associated vertex
algebras are generated by the generating functions. As the
starting point of this paper, for each $m\in \Z$ we form a
generating function
$$E(m,x)=\sum_{n\in \Z}E_{m,m+n}x^{-n}.$$
Then the main defining relations of $\widetilde{\gl}_{\infty}$
can be written as
$$[E(m,x_{1}),E(n,x_{2})]=\left(\frac{x_{1}}{x_{2}}\right)^{m-n}
\left(E(m,x_{2})-E(n,x_{1})+f(m,n){\bf k}\right)$$ for $m,n\in \Z$,
where $f: \Z^{2}\rightarrow \{-1,0,1\}$ is a function determined by
the $2$-cocycle $\psi$ (see Section 2).
The generating functions $E(m,x)$ $(m\in \Z)$ are {\em not} mutually
local, but for any (suitably defined) restricted
$\widetilde{\gl}_{\infty}$-module $W$,
$E(m,x)$ for $m\in \Z$ form what was called in \cite{li-phi} an $S_{trig}$-local set of fields on
$W$.

In order to associate quantum vertex algebras to certain algebras including quantum affine algebras,
a theory of what were called $\phi$-coordinated modules for a
weak quantum vertex algebra was developed  and
a notion of $\S_{trig}$-locality was introduced in \cite{li-phi}, where it was proved
that every $\S_{trig}$-local subset of fields on a vector space $W$
generates in a certain canonical way a (weak) quantum vertex algebra with $W$ as a natural
$\phi$-coordinated module. Taking $W$ to be a (suitably defined) restricted
$\widetilde{\gl}_{\infty}$-module, one can show that
$E(m,x)$ for $m\in \Z$ indeed form an $S_{trig}$-local set of fields on
$W$. In view of this, weak quantum vertex
algebras can be associated to
Lie algebra $\widetilde{\gl}_{\infty}$ {\em conceptually}.

To associate quantum vertex algebras to $\widetilde{\gl}_{\infty}$ {\em explicitly}, we introduce an
infinite dimensional Lie algebra $\widetilde{\gl}_{\infty}^{e}$,
which has generators $B(m,n)$ with $m,n\in \Z$, subject to relations
$$[B(m,x_{1}),B(n,x_{2})]=e^{(m-n)(x_{1}-x_{2})}
\left(B(m,x_{2})-B(n,x_{1})+f(m,n){\bf k}\right),$$ where ${\bf k}$
is a nonzero central element and $B(m,x)=\sum_{n\in
\Z}B(m,n)x^{-n-1}$. For any complex number $\ell$, we construct a
universal ``vacuum module''
$V_{\widetilde{\gl}_{\infty}^{e}}(\ell,0)$ of level $\ell$ for Lie algebra $\widetilde{\gl}_{\infty}^{e}$.
Then by using a result of \cite{li-qva2} we
show that there exists a canonical structure of a quantum vertex
algebra on $V_{\widetilde{\gl}_{\infty}^{e}}(\ell,0)$ and that a restricted
$\widetilde{\gl}_{\infty}^{e}$-module structure of level $\ell$ on a
vector space $W$ exactly amounts to a
$V_{\widetilde{\gl}_{\infty}^{e}}(\ell,0)$-module structure on $W$.
Furthermore, by using \cite{li-phi} we show that a restricted
$\widetilde{\gl}_{\infty}$-module structure of level $\ell$ on a
vector space $W$ exactly amounts to a $\phi$-coordinated
$V_{\widetilde{\gl}_{\infty}^{e}}(\ell,0)$-module structure on $W$.

In literature, there have been many interesting and important
studies on Lie algebras $\gl_{\infty}$ and
$\widetilde{\gl}_{\infty}$, including a remarkable relation with
soliton equations, which was discovered and developed by Kyoto
school (cf. \cite{dkm}, \cite{djkm}), and explicit bosonic and fermionic vertex
operator realizations of the basic (level $1$ irreducible) modules
(see \cite{jm}, \cite{kac1}). Closely related to Lie algebra $\gl_{\infty}$,
certain vertex algebras and their representations have been studied
by several authors (see \cite{fkrw}, \cite{kr}, \cite{xu}). In those
studies, the key idea is to use the embedding of certain Lie algebras into the
completion $\overline{\gl_{\infty}}$. The present study, which
directly uses the Lie algebra
$\widetilde{\gl}_{\infty}$ itself, is different from those in nature.

This paper is organized as follows: In Section 2, we review Lie
algebras $\gl_{\infty}$ and $\widetilde{\gl}_{\infty}$, and
we introduce generating functions and reformulate their defining relations
in terms of generating functions. In Section
3, we introduce Lie algebra $\widetilde{\gl}_{\infty}^{e}$ and using
this we construct a family of quantum vertex algebras
$V_{\widetilde{\gl}_{\infty}^{e}}(\ell,0)$ with a complex parameter
$\ell$. In Section 4, we present a canonical connection between
restricted $\widetilde{\gl}_{\infty}$-modules of level $\ell$ and
$\phi$-coordinated
$V_{\widetilde{\gl}_{\infty}^{e}}(\ell,0)$-modules. In Section 5, we
construct a family of irreducible $\widetilde{\gl}_{\infty}^{e}$-modules.

\section{Lie algebra $\widetilde{\gl}_{\infty}$ and their restricted modules}

This is a short preliminary section. In this section,
we recall from \cite{kac1} the Lie algebra $\gl_{\infty}$
and its $1$-dimensional central extension
$\widetilde{\gl}_{\infty}$, and we formulate a notion of restricted module
and introduce generating functions.

We begin with Lie algebra $\gl_{\infty}$.  By definition,
$\gl_{\infty}$ consists of all doubly infinite complex matrices
$(a_{ij})_{i,j\in \Z}$ with only finitely many nonzero entries. A
canonical base of  $\gl_{\infty}$ consists of matrices $E_{i,j}$ for
$i,j\in \Z$, where $E_{i,j}$ denotes the matrix whose only nonzero
entry is the $(i,j)$-entry which is $1$. We have
\begin{eqnarray}
[E_{m,n},E_{r,s}]=\delta_{n,r}E_{m,s}-\delta_{m,s}E_{r,n}
\end{eqnarray}
for $m,n,r,s\in \Z$. Let $\C^{\infty}$ denote an
infinite-dimensional vector space (over $\C$) with a designated base
$\{ v_{n}\;|\; n\in \Z\}$. Then $\gl_{\infty}$ naturally acts on
$\C^{\infty}$ by
$$E_{i,j}v_{k}=\delta_{j,k}v_{i}\ \ \ \mbox{ for }k\in \Z.$$

Define $\deg E_{i,j}=j-i$ for $i,j\in \Z$ to make $\gl_{\infty}$ a
$\Z$-graded Lie algebra, where the degree-$n$ homogeneous subspace
$(\gl_{\infty})_{(n)}$ for $n\in \Z$ is linearly spanned by
$E_{m,m+n}$ for $m\in \Z$. Using this $\Z$-grading one obtains a
triangular decomposition
$$\gl_{\infty}=\gl_{\infty}^{+}\oplus \gl_{\infty}^{0}\oplus
\gl_{\infty}^{-},$$ where $\gl_{\infty}^{\pm}=\sum_{\pm(j-i)>0}\C
E_{i,j}$ and $\gl_{\infty}^{0}=\sum_{i}\C E_{i,i}$, which is a
Cartan subalgebra. All the trace-zero matrices form a subalgebra
$\sl_{\infty}$, which is isomorphic to the affine Kac-Moody Lie
algebra $\g'(A)$ of type $A_{\infty}$ (see \cite{kac1}). We have
$\gl_{\infty}=\sl_{\infty}\oplus \C E_{0,0}$, where $\gl_{\infty}$
extends $\sl_{\infty}$ by $E_{0,0}$ viewed as a derivation.

On Lie algebra $\gl_{\infty}$, there is a $2$-cocycle $\psi$ defined
by
\begin{eqnarray}
&&\psi(E_{i,j},E_{j,i})=1=-\psi(E_{j,i},E_{i,j})
\ \ \ \mbox{if $i\le 0$ and $j\ge 1$},\nonumber\\
 &&\psi(E_{i,j},E_{m,n})=0 \ \ \ \mbox{otherwise}.
\end{eqnarray}
Using this $2$-cocycle, one obtains a $1$-dimensional central
extension of $\gl_{\infty}$, which is denoted by
$\widetilde{\gl}_{\infty}$. That is,
\begin{eqnarray}
\widetilde{\gl}_{\infty}=\gl_{\infty}\oplus \C {\bf k},
\end{eqnarray}
where ${\bf k}$ is a nonzero central element and
\begin{eqnarray}\label{ehat-form}
[E_{m,n},E_{r,s}]=\delta_{n,r}E_{m,s}-\delta_{m,s}E_{r,n}+\psi(E_{m,n},E_{r,s}){\bf
k}
\end{eqnarray}
for $m,n,r,s\in \Z$.

For convenience, we define a function $f: \Z^{2}\rightarrow
\{-1,0,1\}$ by
\begin{eqnarray}\label{edef-f}
&&f(m,n)=1=-f(n,m)\ \ \ \mbox{if $m\le 0$ and $n\ge 1$},\nonumber\\
 &&f(m,n)=0 \ \ \ \mbox{otherwise}.
\end{eqnarray}
In terms of function $f$ we have
\begin{eqnarray}
\psi(E_{i,j},E_{m,n})=\delta_{i,n}\delta_{j,m}f(i,j)
\ \ \ \mbox{ for }i,j,m,n\in \Z.
\end{eqnarray}

For $m\in \Z$, set
\begin{eqnarray}
E(m,x)=\sum_{n\in \Z}E_{m,m+n}x^{-n}.
\end{eqnarray}

\bl{lgenerating-function} The relations (\ref{ehat-form}) can be
equivalently written as
\begin{eqnarray}\label{edefining-relations}
[E(m,x_{1}),E(n,x_{2})]
=\left(\frac{x_{1}}{x_{2}}\right)^{m-n}\left(E(m,x_{2})-E(n,x_{1})+f(m,n){\bf
k}\right)
\end{eqnarray}
for $m,n\in \Z$. \el

\begin{proof} It is straightforward:
\begin{eqnarray*}
&&[E(m,x_{1}),E(n,x_{2})]\\
& =&\sum_{r,s\in
\Z}\left([E_{m,m+r},E_{n,n+s}]+\delta_{m+r,n}\delta_{m,n+s}f(m,m+r){\bf k}\right)x_{1}^{-r}x_{2}^{-s}\\
&=&\sum_{r,s\in
\Z}\left(\delta_{m+r,n}E_{m,n+s}-\delta_{m,n+s}E_{n,m+r}\right)x_{1}^{-r}x_{2}^{-s}
+f(m,n)\left(\frac{x_{1}}{x_{2}}\right)^{m-n}{\bf k}\\
&=&\left(\frac{x_{1}}{x_{2}}\right)^{m-n}\left(E(m,x_{2})-E(n,x_{1})+f(m,n){\bf
k}\right)
\end{eqnarray*}
for $m,n\in \Z$.
\end{proof}

As an immediate consequence of Lemma \ref{lgenerating-function} we
have:

\bc{cf-function} The function $f: \Z^{2}\rightarrow \{-1,0,1\}$
defined in (\ref{edef-f}) satisfies relations
\begin{eqnarray}
f(m,n)=-f(n,m),\ \ \ \ \ f(m,n)+f(n,r)=f(m,r)
\end{eqnarray}
for $m,n,r\in \Z$. \ec

We formulate the following standard notions:

\bd{dlevel-restricted} {\em A $\widetilde{\gl}_{\infty}$-module $W$
is said to be of {\em level} $\ell\in \C$ if ${\bf k}$ acts on $W$
as scalar $\ell$, and $W$ is said to be {\em restricted} if for
every $m\in \Z$ and for $w\in W$, $E_{m,n}w=0$ for $n$ sufficiently
large.} \ed

For a vector space $W$, set
\begin{eqnarray}
\E(W)=\Hom (W,W((x)))\subset (\End W)[[x,x^{-1}]].
\end{eqnarray}
We see that a $\widetilde{\gl}_{\infty}$-module $W$ is restricted if
and only if $E(m,x)\in \E(W)$ for $m\in \Z$.

It can be readily seen that the natural $\gl_{\infty}$-module
$\C^{\infty}$, which is a $\widetilde{\gl}_{\infty}$-module of level
$0$, is restricted. On the other hand, highest weight
$\widetilde{\gl}_{\infty}$-modules are also restricted modules,
where a highest weight $\widetilde{\gl}_{\infty}$-module with
highest weight $\lambda\in H^{*}$ is a module $W$ with a vector $w$
such that
\begin{eqnarray*}
&&E_{m,m}w=\lambda_{m}w\ \ \ \mbox{ for }m\in \Z,\\
&&E_{m,n}w=0\ \ \ \mbox{ for }m,n\in \Z\ \mbox{ with  }m<n,\\
&&W=U(\widetilde{\gl}_{\infty})w,
\end{eqnarray*}
where $H=\span\{ E_{n,n}\ |\ n\in \Z\}$ and $\lambda_{m}=\lambda(E_{m,m})$.

\section{Lie algebra $\widetilde{\gl}_{\infty}^{e}$ and quantum vertex algebras}
In this section, we introduce an infinite-dimensional Lie algebra
 $\widetilde{\gl}^{e}_{\infty}$, which is intrinsically related to $\widetilde{\gl}_{\infty}$,
and for any complex number $\ell$, we associate a quantum vertex
algebra $V_{\widetilde{\gl}_{\infty}^{e}}(\ell,0)$ to
$\widetilde{\gl}^{e}_{\infty}$ and associate a $V_{\widetilde{\gl}_{\infty}^{e}}(\ell,0)$-module to each
restricted $\widetilde{\gl}^{e}_{\infty}$-module of level $\ell$.

We first recall from \cite{li-qva1} the notion of weak quantum
vertex algebra.

\bd{dweak-qva} {\em A {\em weak quantum vertex algebra} is a vector
space $V$ equipped with a linear map
\begin{eqnarray*}
Y(\cdot,x):&& V\rightarrow \Hom(V,V((x)))\subset (\End
V)[[x,x^{-1}]]\\
&&v\mapsto Y(v,x)=\sum_{n\in \Z}v_{n}x^{-n-1}\ \ (\mbox{where
}v_{n}\in \End V),
\end{eqnarray*}
called the {\em adjoint vertex operator map}, and a vector ${\bf
1}\in V$, called the {\em vacuum vector}, satisfying the following
conditions: For $v\in V$,
$$Y({\bf 1},x)v=v,\ \ Y(v,x){\bf 1}\in V[[x]]\ \ \mbox{and }\
\lim_{x\rightarrow 0}Y(v,x){\bf 1}=v, $$
and for $u,v\in V$, there exist
$$u^{(i)},v^{(i)}\in V,\ \ f_{i}(x)\in \C((x))\ \ \mbox{ for
}i=1,\dots,r$$ such that
\begin{eqnarray}
&&x_{0}^{-1}\delta\left(\frac{x_{1}-x_{2}}{x_{0}}\right)Y(u,x_{1})Y(v,x_{2})\nonumber\\
&&\hspace{1cm}-x_{0}^{-1}\delta\left(\frac{x_{2}-x_{1}}{-x_{0}}\right)
\sum_{i=1}^{r}f_{i}(x_{2}-x_{1})Y(v^{(i)},x_{2})Y(u^{(i)},x_{1})\nonumber\\
&=&x_{2}^{-1}\delta\left(\frac{x_{1}-x_{0}}{x_{2}}\right)Y(Y(u,x_{0})v,x_{2}).
\end{eqnarray}} \ed

For a weak quantum vertex algebra $V$, following \cite{ek} let
$$Y(x): V\otimes V\rightarrow V((x))$$
be the canonical linear map associated to the vertex operator map
$Y(\cdot,x)$.

A {\em rational quantum Yang-Baxter operator} on a vector space $U$
is a linear map
$$\S(x):\ U\otimes U\rightarrow U\otimes U\otimes \C((x)),$$
satisfying
$$\S^{12}(x)\S^{13}(x+z)\S^{23}(z)=\S^{23}(z)\S^{13}(x+z)\S^{12}(x)$$
(the {\em quantum Yang-Baxter equation}), where for $1\le i<j\le 3$,
$$\S^{ij}(x):V\otimes V\otimes V\rightarrow V\otimes V\otimes
V\otimes \C((x))$$ denotes the canonical extension of $\S(x)$.  It
is said to be {\em unitary} if
$$\S(x)\S^{21}(-x)=1,$$
where $\S^{21}(x)=\sigma \S(x)\sigma$ with $\sigma$ denoting the
flip operator on $U\otimes U$.

\bd{dqva} {\em A {\em quantum vertex algebra} is a weak quantum
vertex algebra $V$ equipped with a unitary rational quantum
Yang-Baxter operator $\S(x)$ on $V$, satisfying the following
conditions:
\begin{eqnarray}
&&\S(x)({\bf 1}\otimes v)={\bf 1}\otimes v\ \ \ \mbox{ for }v\in
V,\label{esvacuum}\\
&&[\D\otimes 1, \S(x)]=-\frac{d}{dx}\S(x),\label{d1s}\\
&&Y(u,x)v=e^{x\D}Y(-x)\S(-x)(v\otimes u)\ \ \mbox{ for }u,v\in V,\\
&&\S(x_{1})(Y(x_{2})\otimes 1)=(Y(x_{2})\otimes
1)\S^{23}(x_{1})\S^{13}(x_{1}+x_{2}),\label{sy1}
\end{eqnarray}
where $\D$ is the linear operator on $V$ defined by $\D (v)=v_{-2}{\bf 1}$ for $v\in V$.
We denote a quantum vertex algebra by a pair $(V,\S)$.} \ed

Note that this very notion is a slight modification of the same
named notion in \cite{li-qva1} and \cite{li-qva2} with extra axioms
(\ref{esvacuum}) and (\ref{sy1}).

The following notion was due to Etingof and Kazhdan (see \cite{ek}):

\bd{dkl-nondeg} {\em A weak quantum vertex algebra $V$ is said to be
{\em non-degenerate} if for every positive integer $n$, the linear
map
$$Z_{n}: V^{\otimes n}\otimes \C((x_{1}))\cdots ((x_{n}))\rightarrow
V((x_{1}))\cdots ((x_{n})),$$ defined by
$$Z_{n}(v^{(1)}\otimes \cdots\otimes v^{(n)}\otimes f)
=fY(v^{(1)},x_{1})\cdots Y(v^{(n)},x_{n}){\bf 1}$$ for
$v^{(1)},\dots, v^{(n)}\in V,\; f\in \C((x_{1}))\cdots ((x_{n}))$,
is injective.} \ed

The following is a reformulation of a result in \cite{li-qva1} (cf. \cite{ek}):

\bp{pnon-degenerate} Let $V$ be a weak quantum vertex algebra.
Assume that $V$ is non-degenerate. Then there exists a linear map
$\S(x): V\otimes V\rightarrow V\otimes V\otimes \C((x))$, which is
uniquely determined by
\begin{eqnarray*}
Y(u,x)v=e^{x\D}Y(-x)\S(-x)(v\otimes u) \ \ \ \mbox{for }u,v\in V.
\end{eqnarray*}
Furthermore, $(V,\S)$ carries the structure of a quantum vertex
algebra and the following relation holds
\begin{eqnarray}
[1\otimes \D, \S(x)]=\frac{d}{dx}\S(x).
\end{eqnarray}
 \ep

\br{rnotions} {\em Note that a quantum vertex algebra was defined as
a pair $(V,\S)$. In view of Proposition \ref{pnon-degenerate}, the
term ``a non-degenerate quantum vertex algebra'' without reference
to a quantum Yang-Baxter operator is unambiguous. If a weak quantum vertex algebra
$V$ is of countable dimension over $\C$ and if $V$ as a $V$-module
is irreducible, then by Corollary 3.10 of \cite{li-qva2}, $V$ is
non-degenerate. In view of this, the term ``irreducible quantum
vertex algebra'' is also unambiguous.} \er

\bd{dmodule-wqva} {\em Let $V$ be a weak quantum vertex algebra.
A {\em $V$-module} is a
vector space $W$ equipped with a linear map
$$Y_{W}(\cdot,x): V\rightarrow \Hom (W,W((x)))\subset (\End
W)[[x,x^{-1}]]$$ satisfying the conditions that
$$Y_{W}({\bf
1},x)=1_{W}\ \ (\mbox{where $1_{W}$ denotes the identity operator on
}W)$$ and that for $u,v\in V,\; w\in W$, there exists a nonnegative
integer $l$ such that
\begin{eqnarray}
(x_{0}+x_{2})^{l}Y_{W}(u,x_{0}+x_{2})Y_{W}(v,x_{2})w
=(x_{0}+x_{2})^{l}Y_{W}(Y(u,x_{0})v,x_{2})w.
\end{eqnarray}}
\ed

As we need, we here briefly review the general construction of weak
quantum vertex algebras and their modules from \cite{li-qva1}. Let
$W$ be a vector space. Recall that
$$\E(W)=\Hom
(W,W((x)))\subset (\End W)[[x,x^{-1}]].$$ A subset $U$ of $\E(W)$ is
said to be {\em $\S$-local} if for any $u(x),v(x)\in U$, there exist
$$u^{(i)}(x),v^{(i)}(x)\in U,\ \ f_{i}(x)\in \C((x))
\ \ (i=1,\dots,r)$$ (finitely many) such that
\begin{eqnarray}\label{eslocal-condition}
(x_{1}-x_{2})^{k}u(x_{1})v(x_{2})
=(x_{1}-x_{2})^{k}\sum_{i=1}^{r}f_{i}(x_{2}-x_{1})v^{(i)}(x_{2})u^{(i)}(x_{1})
\end{eqnarray}
for some nonnegative integer $k$.

Notice that the relation (\ref{eslocal-condition}) implies
\begin{eqnarray}\label{ecompatibility-con}
(x_{1}-x_{2})^{k}u(x_{1})v(x_{2})\in \Hom (W,W((x_{1},x_{2}))).
\end{eqnarray}
Now, let $u(x),v(x)\in \E(W)$. Assume that there exists a
nonnegative integer $k$ such that (\ref{ecompatibility-con}) holds.
Define $u(x)_{n}v(x)\in \E(W)$ for $n\in \Z$ in terms of generating
function $$Y_{\E}(u(x),x_{0})v(x)=\sum_{n\in \Z}u(x)_{n}v(x)
x_{0}^{-n-1}$$ by
\begin{eqnarray}
Y_{\E}(u(x),x_{0})v(x)=x_{0}^{-k}\left((x_{1}-x)^{k}u(x_{1})v(x)\right)|_{x_{1}=x+x_{0}}.
\end{eqnarray}
(It was proved that the expression on the right hand side is
independent of the choice of $k$.)  If $u(x),v(x)$ are {}from an
$\S$-local subset $U$, assuming relation (\ref{eslocal-condition})
we have
\begin{eqnarray}
&&u(x)_{n}v(x)\nonumber\\
&=&\Res_{x_{1}}\left((x_{1}-x)^{n}u(x_{1})v(x)-(-x+x_{1})^{n}
\sum_{i=1}^{r}f_{i}(x-x_{1})v^{(i)}(x)u^{(i)}(x_{1})\right).\ \ \ \
\ \ \ \
\end{eqnarray}
The following result was obtained in \cite{li-qva1}:

\bt{tslocal-wqva} Every $\S$-local subset $U$ of $\E(W)$ generates a
weak quantum vertex algebra $\<U\>$ with $W$ as a faithful module
where $Y_{W}(\alpha(x),z)=\alpha(z)$ for $\alpha(x)\in \<U\>$. \et

Next, as one of the main ingredients we introduce a new
infinite-dimensional Lie algebra.

\bp{pliealgebra} Let $E$ be a vector space with basis $\{ e_{m}\ |\
m\in \Z\}$. Set
\begin{eqnarray*}
\widetilde{\gl}_{\infty}^{e}=E\otimes \C((t))\oplus \C {\bf k},
\end{eqnarray*}
where ${\bf k}$ is a symbol. Then there exists a Lie algebra structure on
$\widetilde{\gl}_{\infty}^{e}$ such that
\begin{eqnarray}\label{eBrelation}
&&[{\bf k},\widetilde{\gl}_{\infty}^{e}]=0,\nonumber\\
 &&[B_{t}(m,x_{1}),B_{t}(n,x_{2})]
=e^{(m-n)(x_{1}-x_{2})}\left(B_{t}(m,x_{2})-B_{t}(n,x_{1})+f(m,n){\bf
k}\right)\ \ \ \ \ \ \
\end{eqnarray}
for $m,n\in \Z$, where
\begin{eqnarray}
B_{t}(m,x)=\sum_{n\in \Z}(e_{m}\otimes t^{n})x^{-n-1}.
\end{eqnarray}
\ep

\begin{proof} For $m\in \Z$, we have
$$B_{t}(m,x)=\sum_{n\in \Z}(e_{m}\otimes t^{n})x^{-n-1}
=e_{m}\otimes x^{-1}\delta\left(\frac{t}{x}\right).$$ Furthermore,
for $g(t)\in \C((t))$ we have
$$g(x)B_{t}(m,x)=e_{m}\otimes g(t)x^{-1}\delta\left(\frac{t}{x}\right)
\in \widetilde{\gl}_{\infty}^{e}[[x,x^{-1}]]$$ and
\begin{eqnarray}
e_{m}\otimes g(t)=\Res_{x}g(x)B_{t}(m,x).
\end{eqnarray}

Define a bilinear operation $[\cdot,\cdot]$ on
$\widetilde{\gl}_{\infty}^{e}$ by
$$\left[{\bf
k},\widetilde{\gl}_{\infty}^{e}\right]=0=\left[\widetilde{\gl}_{\infty}^{e},{\bf
k}\right],$$
\begin{eqnarray}\label{edef-bracket}
&&\left[e_{m}\otimes g(t),e_{n}\otimes
h(t)\right]\nonumber\\
&=&\Res_{x_{1}}\Res_{x_{2}}g(x_{1})h(x_{2})
e^{(m-n)(x_{1}-x_{2})}\left(B_{t}(m,x_{2})-B_{t}(n,x_{1})+f(m,n){\bf
k}\right)\ \ \ \ \ \ \ \
\end{eqnarray}
for $m,n\in \Z,\; g(t),h(t)\in \C((t))$. Indeed, for any fixed
integers $m,n$, it is $\C$-linear in both $g(t)$ and $h(t)$.
It can be readily seen that skew symmetry holds as the function
$f(m,n)$ on $\Z\times \Z$ is skew symmetric. We next show that
Jacobi identity also holds.

{}From definition, (\ref{eBrelation}) holds for $m,n\in \Z$.
Furthermore, by the standard formal variable convention we have
\begin{eqnarray}
&&\left[g(x_{1})B_{t}(m,x_{1}),h(x_{2})B_{t}(n,x_{2})\right]\nonumber\\
&=&g(x_{1})h(x_{2})e^{(m-n)(x_{1}-x_{2})}\left(B_{t}(m,x_{2})-B_{t}(n,x_{1})+f(m,n){\bf
k}\right)
\end{eqnarray}
for $g(t),h(t)\in \C((t))$.

Let $m,n,r\in \Z$. We have
\begin{eqnarray*}
&&\left[[B_{t}(m,x_{1}),B_{t}(n,x_{2})],B_{t}(r,x_{3})\right]\\
&=&e^{(m-n)(x_{1}-x_{2})}\left[B_{t}(m,x_{2})-B_{t}(n,x_{1})+f(m,n){\bf
k},
B_{t}(r,x_{3})\right]\\
&=&e^{(m-n)(x_{1}-x_{2})}
e^{(m-r)(x_{2}-x_{3})}\left(B_{t}(m,x_{3})-B_{t}(r,x_{2})+f(m,r){\bf k}\right)\\
&&-e^{(m-n)(x_{1}-x_{2})}
e^{(n-r)(x_{1}-x_{3})}\left(B_{t}(n,x_{3})-B_{t}(r,x_{1})+f(n,r){\bf
k}\right)\\
&=&e^{(m-n)x_{1}+(n-r)x_{2}+(r-m)x_{3}}
\left(B_{t}(m,x_{3})-B_{t}(r,x_{2})+f(m,r){\bf k}\right)\\
&&-e^{(m-r)x_{1}+(n-m)x_{2}+(r-n)x_{3}}
\left(B_{t}(n,x_{3})-B_{t}(r,x_{1})+f(n,r){\bf k}\right),
\end{eqnarray*}
\begin{eqnarray*}
&&\left[B_{t}(m,x_{1}),[B_{t}(n,x_{2}),B_{t}(r,x_{3})]\right]\\
&=& e^{(n-r)(x_{2}-x_{3})}
e^{(m-n)(x_{1}-x_{3})}\left(B_{t}(m,x_{3})-B_{t}(n,x_{1})+f(m,n){\bf
k}\right)\\
&&-e^{(n-r)(x_{2}-x_{3})}
e^{(m-r)(x_{1}-x_{2})}\left(B_{t}(m,x_{2})-B_{t}(r,x_{1})+f(m,r){\bf
k}\right)\\
&=&e^{(m-n)x_{1}+(n-r)x_{2}+(r-m)x_{3}}\left(B_{t}(m,x_{3})-B_{t}(n,x_{1})+f(m,n){\bf
k}\right)\\
&&-e^{(m-r)x_{1}+(n-m)x_{2}+(r-n)x_{3}}
\left(B_{t}(m,x_{2})-B_{t}(r,x_{1})+f(m,r){\bf k}\right),
\end{eqnarray*}
\begin{eqnarray*}
&&\left[B_{t}(n,x_{2}),[B_{t}(m,x_{1}),B_{t}(r,x_{3})]\right]\\
&=&e^{(m-r)x_{1}+(n-m)x_{2}+(r-n)x_{3}}
\left(B_{t}(n,x_{3})-B_{t}(m,x_{2})+f(n,m){\bf
k}\right)\\
&&-e^{(m-n)x_{1}+(n-r)x_{2}+(r-m)x_{3}}
\left(B_{t}(n,x_{1})-B_{t}(r,x_{2})+f(n,r){\bf k}\right).
\end{eqnarray*}
Recall from Corollary \ref{cf-function} that
\begin{eqnarray*}
f(m,n)+f(n,r)=f(m,r),\ \ \ f(m,r)+f(n,m)=f(n,r).
\end{eqnarray*}
(Note that the second relation follows from the first one as
$f(m,n)=-f(n,m)$.) Then we see that the Jacobi identity holds.
\end{proof}

\br{rphi-mod} {\em For $n\in \Z$, set $U_{n}=E\otimes t^{n}\C[[t]]$.
This gives a descending sequence $\{ U_{n}\}_{n\in \Z}$ with
$\cap_{n\in \Z}U_{n}=0$. Equip vector space
$\widetilde{\gl}_{\infty}^{e}$ with the topological vector space
structure associated to $\{ U_{n}\}_{n\in \Z}$. Let $u\in E,\;
g(t)\in \C((t))$. From definition (see (\ref{edef-bracket})) we have
\begin{eqnarray}\label{elie-continous}
[u\otimes g(t),U_{n}]\subset U_{n}\ \ \ \mbox{ for  }n\ge 0,
\end{eqnarray}
which implies that the Lie bracket is continuous. Thus
$\widetilde{\gl}_{\infty}^{e}$ is a topological Lie algebra where
the elements ${\bf k}$ and $e_{m}\otimes t^{r}$ for $m,r\in \Z$ form
a topological basis. } \er

\br{rfiltered} {\em  For $n\in \Z$, set
\begin{eqnarray}
\widetilde{\gl}_{\infty}^{e}[n]=\begin{cases}E\otimes t^{n}\C[[t]]\
\
\mbox{ if }n\ge 1\\
E\otimes t^{n}\C[[t]]+\C {\bf k}\ \ \mbox{ if }n\le 0.
\end{cases}
\end{eqnarray}
This defines a descending filtration of
$\widetilde{\gl}_{\infty}^{e}$. Using (\ref{edef-bracket}) one can
show that
$$\left[\widetilde{\gl}_{\infty}^{e}[m],\widetilde{\gl}_{\infty}^{e}[n]\right]\subset
\widetilde{\gl}_{\infty}^{e}[m+n]$$ for $m,n\in \Z$. Thus,
$\widetilde{\gl}_{\infty}^{e}$ equipped with this filtration becomes
a filtered Lie algebra. (But, Lie algebra
$\widetilde{\gl}_{\infty}^{e}$ is not $\Z$-graded in the obvious
way.)} \er

For $m,r\in \Z$, let $B(m,r)$ denote the operator corresponding to
$e_{m}\otimes t^{r}$ on a $\widetilde{\gl}_{\infty}^{e}$-module. For $m\in \Z$, set
\begin{eqnarray}
B(m,x)=\sum_{n\in \Z}B(m,n)x^{-n-1}.
\end{eqnarray}
We define a notion of restricted
$\widetilde{\gl}_{\infty}^{e}$-module in the usual way. That is, a
$\widetilde{\gl}_{\infty}^{e}$-module $W$ is said to be {\em
restricted} if for any $m\in \Z,\ w\in W$, $B(m,n)w=0$ for $n$
sufficiently large, or equivalently $B(m,x)\in \E(W)$ for $m\in \Z$.
We also assume continuity.

Set
$${\mathcal{B}}^{+}=E\otimes \C[[t]]\subset \widetilde{\gl}_{\infty}^{e}.$$
Let $m,n\in \Z,\; g(t), h(t)\in \C[[t]]$. Noticing that
$e^{(m-n)(x_{1}-x_{2})}\in \C[[x_{1},x_{2}]]$, from
(\ref{edef-bracket}) we get
\begin{eqnarray*}
[e_{m}\otimes g(t),e_{n}\otimes h(t)]=0.
\end{eqnarray*}
Thus, ${\mathcal{B}}^{+}$ is an abelian subalgebra of
$\widetilde{\gl}_{\infty}^{e}$. Denote by ${\mathcal{B}}^{-}$ the
subspace of $\widetilde{\gl}_{\infty}^{e}$, linearly spanned by $e_{m}\otimes
t^{r}$ for $m,r\in \Z$ with $r<0$. We have
$$\widetilde{\gl}_{\infty}^{e}
={\mathcal{B}}^{+}\oplus \C {\bf k}\oplus {\mathcal{B}}^{-}$$ as a
vector space. (Note that ${\mathcal{B}}^{-}$ is not a subalgebra.)

Let $\ell$ be any complex number. Letting ${\mathcal{B}}^{+}$ act on $\C$
trivially and letting ${\bf k}$ act as scalar $\ell$, we obtain a
$({\mathcal{B}}^{+}\oplus \C {\bf k})$-module denoted by
$\C_{\ell}$. Then form an induced module
\begin{eqnarray}
V_{\widetilde{\gl}_{\infty}^{e}}(\ell,0)=U(\widetilde{\gl}_{\infty}^{e})\otimes
_{U({\mathcal{B}}^{+}\oplus \C {\bf k})}\C_{\ell}.
\end{eqnarray}
In view of the P-B-W theorem,
$V_{\widetilde{\gl}_{\infty}^{e}}(\ell,0)$ as a vector space is
isomorphic to $S({\mathcal{B}}^{-})$ (the symmetric algebra over
${\mathcal{B}}^{-}$). Set ${\bf 1}=1\otimes 1\in
V_{\widetilde{\gl}_{\infty}^{e}}(\ell,0)$, and then set
\begin{eqnarray}
b^{(m)}=B(m,-1){\bf 1}\in V_{\widetilde{\gl}_{\infty}^{e}}(\ell,0) \
\ \mbox{ for }m\in \Z.
\end{eqnarray}

Note that $B(m,k){\bf 1}=0$ for $m\in \Z,\; k\in \N$. Let $m,n\in
\Z,\; k\in \N$. {}From the Lie bracket relation (\ref{eBrelation}) we
get
\begin{eqnarray}
[B(m,k), B(n,x_{2})]=-e^{(n-m)x_{2}}\sum_{j\ge
0}\frac{(m-n)^{j}}{j!}B(n,k+j).
\end{eqnarray}
It follows from this relation and induction that $B(m,k)=0$ on
$V_{\widetilde{\gl}_{\infty}^{e}}(\ell,0)$ for $m\in \Z,\; k\ge 0$.
In particular, $V_{\widetilde{\gl}_{\infty}^{e}}(\ell,0)$ is a
restricted
 $\widetilde{\gl}_{\infty}^{e}$-module.

As one of the main results in this section, we have:

\bt{tqva} For every complex number $\ell$, there exists a weak
quantum vertex algebra structure on
$V_{\widetilde{\gl}_{\infty}^{e}}(\ell,0)$, which is uniquely
determined by the conditions that ${\bf 1}$ is the vacuum vector and
that
$$Y(b^{(m)},x)=B(m,x)\ \ \mbox{
for }m\in \Z.$$ Furthermore,
$V_{\widetilde{\gl}_{\infty}^{e}}(\ell,0)$ is a non-degenerate
quantum vertex algebra and the following relations hold:
\begin{eqnarray}\label{esingular}
b^{(m)}_{k}b^{(n)}=0\ \ \ \mbox{ for }m,n\in \Z, \ k\in \N
\end{eqnarray}
 and
\begin{eqnarray}\label{evertex-operator-relations}
&&Y(b^{(m)},x_{1})Y(b^{(n)},x_{2})-Y(b^{(n)},x_{2})Y(b^{(m)},x_{1})\nonumber\\
&=&e^{(m-n)(x_{1}-x_{2})}\left(Y(b^{(m)},x_{2})-Y(b^{(n)},x_{1})+f(m,n)\ell\right)
\end{eqnarray}
for $m,n\in \Z$. \et

\begin{proof} Let $W$ be any restricted $\widetilde{\gl}_{\infty}^{e}$-module of
level $\ell$. Set $$U_{W}=\{ B(m,x)\;|\; m\in \Z\}\cup \{
1_{W}\}\subset \E(W).$$ Writing the defining relation
(\ref{eBrelation}) as
\begin{eqnarray}\label{ereformulation1}
&&B(m,x_{1})B(n,x_{2})=B(n,x_{2})B(m,x_{1})\nonumber\\
&&\ \ \ \ \
+e^{(m-n)(x_{1}-x_{2})}\left(B(m,x_{2})-B(n,x_{1})+f(m,n)\ell\right),
\end{eqnarray}
we see that $U_{W}$ is an $\S$-local subset of $\E(W)$. Then $U_{W}$
generates a weak quantum vertex algebra $\<U_{W}\>$ inside $\E(W)$
with $W$ as a canonical module. With (\ref{ereformulation1}), by
Proposition 3.13 of \cite{li-qva1} (cf. Proposition 2.12,
\cite{li-qva2}), we have
\begin{eqnarray*}
&&Y_{\E}(B(m,x),x_{1})Y_{\E}(B(n,x),x_{2})
=Y_{\E}(B(n,x),x_{2})Y_{\E}(B(m,x)x_{1})\\
&&\ \ \ \ +
e^{(m-n)(x_{1}-x_{2})}\left(Y_{\E}(B(m,x),x_{2})-Y_{\E}(B(n,x),x_{1})+f(m,n)\ell\right)
\end{eqnarray*}
for $m,n\in \Z$.  This shows that $\<U_{W}\>$ is a
$\widetilde{\gl}_{\infty}^{e}$-module of level $\ell$ with $B(m,z)$
acting as $Y_{\E}(B(m,x),z)$ for $m\in \Z$. Furthermore, we have
$$B(m,x)_{k}1_{W}=0\ \ \ \mbox{ for }m\in \Z,\; k\in \N.$$
It follows from the construction of
$V_{\widetilde{\gl}_{\infty}^{e}}(\ell,0)$ that there exists a
$\widetilde{\gl}_{\infty}^{e}$-module homomorphism $\theta$ from
$V_{\widetilde{\gl}_{\infty}^{e}}(\ell,0)$ to $\<U_{W}\>$ with
$\theta({\bf 1})=1_{W}$.

Take $W=V_{\widetilde{\gl}_{\infty}^{e}}(\ell,0)$. Then it follows
from Theorem 2.9 of \cite{li-qva2} that there is a weak quantum
vertex algebra structure on
$V_{\widetilde{\gl}_{\infty}^{e}}(\ell,0)$ with all the desired
properties.

As for non-degeneracy, we shall use a result of \cite{li-qva2}. For
$n\ge 0$, let $V[n]$ denote the subspace of
$V_{\widetilde{\gl}_{\infty}^{e}}(\ell,0)$ linearly spanned by vectors
$$b^{(m_{1})}_{k_{1}}\cdots b^{(m_{r})}_{k_{r}}{\bf 1}$$
where $0\le r\le n,\ m_{i},k_{i}\in \Z$ for $i=1,\dots,r$. (Note
that $b^{(m)}_{k}=B(m,k)$.) Since $\{ b^{(m)}\ |\ m\in \Z\}$
generates $V_{\widetilde{\gl}_{\infty}^{e}}(\ell,0)$ as a weak
quantum vertex algebra, this defines an ascending filtration of
$V_{\widetilde{\gl}_{\infty}^{e}}(\ell,0)$. Denote the associated
graded weak quantum vertex algebra by $K$. By Proposition 3.15 of
\cite{li-qva2}, $K$ is a commutative vertex algebra. Furthermore, it
follows {}from the P-B-W basis that $K=S({\mathcal{B}}^{-})$. Then,
by Theorem 3.19 of \cite{li-qva2},
$V_{\widetilde{\gl}_{\infty}^{e}}(\ell,0)$ is non-degenerate.

Since $Y(b^{(m)},x)=B(m,x)$ for $m\in \Z$, the last assertion is
clear. For $m,n\in \Z$, from the last assertion we have
\begin{eqnarray*}
&&Y(b^{(m)},x_{1})Y(b^{(n)},x_{2})=Y(b^{(n)},x_{2})Y(b^{(m)},x_{1})\\
&&+e^{(m-n)(x_{1}-x_{2})}\left( Y(b^{(m)},x_{2})Y({\bf
1},x_{1})-Y({\bf 1},x_{2})Y(b^{(n)},x_{1}) +\ell f(m,n)Y({\bf
1},x_{2})Y({\bf 1},x_{1})\right),
\end{eqnarray*}
noticing that $Y({\bf 1},x)=1$. Set
\begin{eqnarray*}
&&P=Y(b^{(n)},x_{2})Y(b^{(m)},x_{1})\\
&&+e^{(m-n)(x_{1}-x_{2})}\left( Y(b^{(m)},x_{2})Y({\bf
1},x_{1})-Y({\bf 1},x_{2})Y(b^{(n)},x_{1}) +\ell f(m,n)Y({\bf
1},x_{2})Y({\bf 1},x_{1})\right).\ \ \ \ \
\end{eqnarray*}
For $k\ge 0$, from the $\S$-Jacobi identity for weak quantum vertex
algebras we have
\begin{eqnarray*}
Y(b^{(m)}_{k}b^{(n)},x_{2})
=\Res_{x_{1}}(x_{1}-x_{2})^{k}\left(Y(b^{(m)},x_{1})Y(b^{(n)},x_{2})-P\right)
=0.
\end{eqnarray*}
It then follows that $b^{(m)}_{k}b^{(n)}=0$ for $m,n\in \Z,\; k\in
\N$.
\end{proof}

Furthermore, we have:

\bt{tqva-module} For any restricted
$\widetilde{\gl}_{\infty}^{e}$-module $W$ of level $\ell$, there
exists a $V_{\widetilde{\gl}_{\infty}^{e}}(\ell,0)$-module structure
$Y_{W}$, which is uniquely determined by $$Y_{W}(b^{(m)},x)=B(m,x)\
\ \ \mbox{ for }m\in \Z.$$ On the other hand, for any
$V_{\widetilde{\gl}_{\infty}^{e}}(\ell,0)$-module $(W,Y_{W})$, $W$
is a restricted $\widetilde{\gl}_{\infty}^{e}$-module of level
$\ell$ with
$$B(m,x)=Y_{W}(b^{(m)},x)\ \ \ \mbox{ for }m\in \Z.$$
\et

\begin{proof} Set $U_{W}=\{ B(m,x)\;|\; m\in \Z\}\cup \{ 1_{W}\}\subset \E(W)$.
{}From the proof of Theorem \ref{tqva}, we have a weak quantum
vertex algebra $\<U_{W}\>$ generated by $U_{W}$ and $W$ is a
$\<U_{W}\>$-module with $Y_{W}(\alpha(x),z)=\alpha(z)$ for
$\alpha(x)\in \<U_{W}\>$. Furthermore, $\<U_{W}\>$ is a
$\widetilde{\gl}_{\infty}^{e}$-module of level $\ell$ with
$B(m,z)=Y_{\E}(B(m,x),z)$ for $m\in \Z$ and there exists a
$\widetilde{\gl}_{\infty}^{e}$-module homomorphism $\theta$ from
$V_{\widetilde{\gl}_{\infty}^{e}}(\ell,0)$ to $\<U_{W}\>$ with
$\theta({\bf 1})=1_{W}$.  We have
$$\theta(b^{(m)})=\theta(B(m,-1){\bf 1})=B(m,x)_{-1}1_{W}=B(m,x)$$
and
$$\theta(Y(b^{(m)},z)v)=\theta(B(m,z)v)=Y_{\E}(B(m,x),z)\theta(v)=Y_{\E}(\theta(b^{(m)}),z)\theta(v)$$
for $m\in \Z,\; v\in V_{\widetilde{\gl}_{\infty}^{e}}(\ell,0)$. As
$\{ b^{(m)}\; |\; m\in \Z\}$ generates
$V_{\widetilde{\gl}_{\infty}^{e}}(\ell,0)$ as a weak quantum vertex
algebra, it follows that $\theta$ is a homomorphism of weak quantum
vertex algebras. Consequently, $W$ becomes a
$V_{\widetilde{\gl}_{\infty}^{e}}(\ell,0)$-module with
$$Y_{W}(v,x)=\theta(v)(x)\in \<U_{W}\>\subset \E(W)$$
for $v\in V_{\widetilde{\gl}_{\infty}^{e}}(\ell,0)$. In particular,
we have $Y_{W}(b^{(m)},x)=\theta(b^{(m)})=B(m,x)$ for $m\in \Z$.

On the other hand, let $(W,Y_{W})$ be a
$V_{\widetilde{\gl}_{\infty}^{e}}(\ell,0)$-module. With
(\ref{evertex-operator-relations}), by Corollary 5.4 of
\cite{li-qva1}, we have
\begin{eqnarray*}
&&Y_{W}(b^{(m)},x_{1})Y_{W}(b^{(n)},x_{2})-Y_{W}(b^{(n)},x_{2})Y_{W}(b^{(m)},x_{1})\nonumber\\
&=&e^{(m-n)(x_{1}-x_{2})}\left(Y_{W}(b^{(m)},x_{2})-Y_{W}(b^{(n)},x_{1})+f(m,n)\ell\right)
\end{eqnarray*}
for $m,n\in \Z$. That is, $W$ becomes a
$\widetilde{\gl}_{\infty}^{e}$-module of level $\ell$ with
$B(m,x)=Y_{W}(b^{(m)},x)$ for $m\in \Z.$ It is restricted as
$Y_{W}(b^{(m)},x)\in \E(W)$ for $m\in \Z$ from definition.
\end{proof}

\section{Lie algebra $\widetilde{\gl}_{\infty}$ and
quantum vertex algebra $V_{\widetilde{\gl}_{\infty}^{e}}(\ell,0)$}

In this section, we relate restricted
$\widetilde{\gl}_{\infty}$-modules of level $\ell$ with the quantum
vertex algebra $V_{\widetilde{\gl}_{\infty}^{e}}(\ell,0)$ in terms
of $\phi$-coordinated modules in the sense of \cite{li-phi}. More
specifically, we show that a level-$\ell$ restricted
$\widetilde{\gl}_{\infty}$-module structure on a vector space $W$
exactly amounts to a $\phi$-coordinated module structure  for the
quantum vertex algebra $V_{\widetilde{\gl}_{\infty}^{e}}(\ell,0)$.

We first recall from \cite{li-phi} some basic notions and results
on $\phi$-coordinated modules for weak quantum vertex algebras.
 Set
$$\phi=\phi(x,z)=xe^{z}\in \C[[x,z]],$$ which is fixed throughout
this section.

\bd{dphi-module} {\em Let $V$ be a weak quantum vertex algebra. A
{\em $\phi$-coordinated $V$-module} is a vector space $W$, equipped
with a linear map
$$Y_{W}(\cdot,x): V\rightarrow \Hom (W,W((x)))\subset (\End
W)[[x,x^{-1}]],$$ satisfying the conditions that $$Y_{W}({\bf
1},x)=1_{W}$$ and that for $u,v\in V$, there exists a nonnegative
integer $k$ such that
\begin{eqnarray}
(x_{1}-x_{2})^{k}Y_{W}(u,x_{1})Y_{W}(v,x_{2})\in \Hom
(W,W((x_{1},x_{2})))
\end{eqnarray}
 and
\begin{eqnarray}
(e^{x_{0}}-1)^{k}Y_{W}(Y(u,x_{0})v,x_{2})
=\left((x_{1}/x_{2}-1)^{k}Y_{W}(u,x_{1})Y_{W}(v,x_{2})\right)|_{x_{1}=x_{2}e^{x_{0}}}.
\end{eqnarray}}
\ed

Let $\C(x)$ denote the field of rational functions in $x$. That is,
$\C(x)$ is the fraction field of the polynomial algebra $\C[x]$. On
the other hand, let $\C_{*}(x)$ denote the fraction field of the
algebra $\C[[x]]$. There exists an algebra map $\iota:
\C_{*}(x)\rightarrow \C((x))$, which is uniquely determined by
$$\iota (q(x))=q(x)\ \ \ \mbox{ for }q(x)\in \C[[x]].$$
For $f(x)\in \C(x)$, $\iota (f(x))$ is simply the formal Laurent
series expansion of $f(x)$ at $x=0$. For $f(x)\in \C(x)$, we have
$$(\iota f)(x_{1}/x_{2})=\iota_{x_{2},x_{1}}f(x_{1}/x_{2}).$$

Let $W$ be a vector space. Recall that
$$\E(W)=\Hom (W,W((x)))\subset (\End W)[[x,x^{-1}]].$$
Let $a(x),b(x)\in \E(W)$. Assume that there exists a nonzero
$p(x)\in \C[x]$ such that
\begin{eqnarray}\label{eprecondition}
p(x_{1}/x_{2})a(x_{1})b(x_{2})\in \Hom (W,W((x_{1},x_{2}))).
\end{eqnarray}
Define $a(x)_{n}^{e}b(x)\in \E(W)$ for $n\in \Z$ in terms of
generating function $$Y_{\E}^{e}(a(x),z)b(x)=\sum_{n\in
\Z}a(x)_{n}^{e}b(x) z^{-n-1}$$ by
\begin{eqnarray}
Y_{\E}^{e}(a(x),z)b(x)
=\iota(1/p(e^{z}))\left(p(x_{1}/x)a(x_{1})b(x)\right)|_{x_{1}=xe^{z}},
\end{eqnarray}
where $p(x)$ is any nonzero polynomial such that
(\ref{eprecondition}) holds. (It was proved in \cite{li-phi} that
$p(e^{z})\ne 0$ in $\C[[z]]$ for any nonzero polynomial $p(x)$.)

\bd{dslocal} {\em A subset $U$ of $\E(W)$ is said to be {\em
$\S_{trig}$-local} if for any $u(x),v(x)\in U$, there exist
$$u^{(i)}(x),v^{(i)}(x)\in U,\ \ f_{i}(x)\in \C(x)
\ \ (i=1,\dots,r)$$ (finitely many) such that
\begin{eqnarray}\label{estrig-local}
(x_{1}-x_{2})^{k}u(x_{1})v(x_{2})
=(x_{1}-x_{2})^{k}\sum_{i=1}^{r}(\iota
f_{i})(x_{1}/x_{2})v^{(i)}(x_{2})u^{(i)}(x_{1})
\end{eqnarray}
for some nonnegative integer $k$.} \ed

Notice that relation (\ref{estrig-local}) implies that
$$(x_{1}/x_{2}-1)^{k}u(x_{1})v(x_{2})\in \Hom(W,W((x_{1},x_{2}))).$$
The following result was obtained in \cite{li-phi}:

\bt{trecall-1} Let $W$ be a vector space and let $U$ be an
$\S_{trig}$-local subset of $\E(W)$. Then $U$ generates a weak
quantum vertex algebra $\<U\>_{e}$ with $W$ as a faithful
$\phi$-coordinated module where $Y_{W}(\alpha(x),z)=\alpha(z)$ for
$\alpha(x)\in \<U\>_{e}$. \et

Now, we are in a position to present the main result of this paper.

\bt{tphi-module} Let $\ell$ be any complex number. For every
restricted $\widetilde{\gl}_{\infty}$-module $W$ of level $\ell$,
there exists a $\phi$-coordinated
$V_{\widetilde{\gl}_{\infty}^{e}}(\ell,0)$-module structure $Y_{W}$
on $W$, which is uniquely determined by
$$Y_{W}(b^{(m)},x)=E(m,x)\ \ \mbox{ for }m\in \Z.$$ On the other
hand, for every $\phi$-coordinated
$V_{\widetilde{\gl}_{\infty}^{e}}(\ell,0)$-module $(W,Y_{W})$, $W$
is a restricted $\widetilde{\gl}_{\infty}$-module of level $\ell$
with
$$E(m,x)=Y_{W}(b^{(m)},x)\ \ \mbox{
for }m\in \Z.$$ \et

\begin{proof} Set
$$U_{W}=\{ E(m,x)\ |\ m\in \Z\}\cup \{ 1_{W}\}\subset \E(W).$$
It can be readily seen from (\ref{edefining-relations}) that $U_{W}$
is an $\S_{trig}$-local subset of $\E(W)$. By Theorem
\ref{trecall-1}, $U_{W}$ generates a weak quantum vertex algebra
$\<U_{W}\>_{e}$, where $1_{W}$ serves as the vacuum vector and the
vertex operator map is denoted by $Y_{\E}^{e}(\cdot,x)$.
Furthermore, $W$ becomes a $\phi$-coordinated module with
$Y_{W}(\alpha(x),z)=\alpha(z)$ for $\alpha(x)\in \<U_{W}\>_{e}$.

With the relations (\ref{edefining-relations}) (with ${\bf
k}=\ell$), by Proposition 5.3 of \cite{li-phi}, we have
\begin{eqnarray*}
&&Y_{\E}^{e}(E(m,x),x_{1})Y_{\E}^{e}(E(n,x),x_{2})
-Y_{\E}^{e}(E(n,x),x_{2})Y_{\E}^{e}(E(m,x),x_{1})\\
&=&e^{(m-n)(x_{1}-x_{2})}\left(Y_{\E}^{e}(E(m,x),x_{2})-Y_{\E}^{e}(E(n,x),x_{1})+f(m,n)\ell
1_{W}\right)
\end{eqnarray*}
for $m,n\in \Z$. This shows that $\<U_{W}\>_{e}$ becomes a
restricted $\widetilde{\gl}_{\infty}^{e}$-module of level $\ell$
with $B(m,z)$ acting as $Y_{\E}^{e}(E(m,x),z)$ for $m\in \Z$.
Furthermore, just as in the proof of Theorem \ref{tqva}, we have
that $\<U_{W}\>_{e}$ as a $\widetilde{\gl}_{\infty}^{e}$-module is
cyclic on $1_{W}$ and
$$E(m,x)^{e}_{k}1_{W}=0\ \ \ \mbox{ for }m\in \Z,\; k\in \N.$$
It follows from the construction of the
$\widetilde{\gl}_{\infty}^{e}$-module
$V_{\widetilde{\gl}_{\infty}^{e}}(\ell,0)$ that there exists a
$\widetilde{\gl}_{\infty}^{e}$-module homomorphism $\theta$ from
$V_{\widetilde{\gl}_{\infty}^{e}}(\ell,0)$ to $\<U_{W}\>_{e}$ with
$\theta({\bf 1})=1_{W}$. For $m\in \Z$, we have
$$\theta(b^{(m)})=\theta(B(m,-1){\bf
1})=B(m,-1)\theta({\bf 1})=E(m,x)_{-1}^{e}1_{W}=E(m,x).$$
Furthermore, we have
$$\theta(Y(b^{(m)},z)v)=\theta(B(m,z)v)
=Y_{\E}^{e}(E(m,x),z)\theta(v)=Y_{\E}^{e}(\theta(b^{(m)}),z)\theta(v)$$
for $m\in \Z,\ v\in V_{\widetilde{\gl}_{\infty}^{e}}(\ell,0)$.
 Since $\{ b^{(m)}\ |\ m\in \Z\}$
generates $V_{\widetilde{\gl}_{\infty}^{e}}(\ell,0)$ as a weak
quantum vertex algebra, it follows that $\theta$ is a homomorphism
of weak quantum vertex algebras. With $W$ as a canonical
$\phi$-coordinated module for $\<U_{W}\>_{e}$, through the
homomorphism $\theta$, $W$ becomes a $\phi$-coordinated
$V_{\widetilde{\gl}_{\infty}^{e}}(\ell,0)$-module, where
$$Y_{W}(b^{(m)},x)=Y_{W}(\theta(b^{(m)}),x)=E(m,x)\ \ \ \mbox{ for
}m\in \Z.$$ The uniqueness follows from the fact that $\{ b^{(m)}\
|\ m\in \Z\}$ generates $V_{\widetilde{\gl}_{\infty}^{e}}(\ell,0)$
as a weak quantum vertex algebra.

On the other hand, let $(W,Y_{W})$ be a $\phi$-coordinated
$V_{\widetilde{\gl}_{\infty}^{e}}(\ell,0)$-module. For $m,n\in \Z$,
recalling (\ref{evertex-operator-relations}), {}from \cite{li-phi}
(Proposition 5.9) we have
\begin{eqnarray*}
&&Y_{W}(b^{(m)},x_{1})Y_{W}(b^{(n)},x_{2})-Y_{W}(b^{(n)},x_{2})Y_{W}(b^{(m)},x_{1})\\
&&\ \ \ \ -\left(\frac{x_{1}}{x_{2}}\right)^{m-n}\left(
Y_{W}(b^{(m)},x_{2})-Y_{W}(b^{(n)},x_{1})+\ell
f(m,n)\right)\\
&=&\sum_{k\ge
0}Y_{W}(b^{(m)}_{k}b^{(n)},x_{2})\frac{1}{k!}\left(x_{2}\frac{\partial}{\partial
x_{2}}\right)^{k}\delta\left(\frac{x_{2}}{x_{1}}\right).
\end{eqnarray*}
Recall that for quantum vertex algebra
$V_{\widetilde{\gl}_{\infty}^{e}}(\ell,0)$, the following relations
hold:
$$b^{(m)}_{k}b^{(n)}=0\ \ \ \mbox{ for }k\ge 0.$$
Then
\begin{eqnarray*}
&&Y_{W}(b^{(m)},x_{1})Y_{W}(b^{(n)},x_{2})
-Y_{W}(b^{(n)},x_{2})Y_{W}(b^{(m)},x_{1})\nonumber\\
&=&\left(\frac{x_{1}}{x_{2}}\right)^{m-n}\left(
Y_{W}(b^{(m)},x_{2})-Y_{W}(b^{(n)},x_{1})+\ell f(m,n)\right).
\end{eqnarray*}
Thus $W$ becomes a $\widetilde{\gl}_{\infty}$-module of level $\ell$
with $E(m,x)$ acting as $Y_{W}(b^{(m)},x)$ for $m\in \Z$, and
clearly it is a restricted module.
\end{proof}

\section{Constructing $\widetilde{\gl}_{\infty}^{e}$-modules}
In this section, we  give a construction
of $\widetilde{\gl}_{\infty}^{e}$-modules {}from
$\widetilde{\gl}_{\infty}$-modules of a certain type, including the
natural module $\C^{\infty}$. We also discuss a
generalized-Verma-module construction of
$\widetilde{\gl}_{\infty}^{e}$-modules.

First, we introduce a special family of restricted
$\widetilde{\gl}_{\infty}$-modules, including the natural module
$\C^{\infty}$ (of level $0$).

\bd{dcfin-2} {\em Let $C_{\rm fin}$ denote the category of
$\widetilde{\gl}_{\infty}$-modules $W$ such that
\begin{eqnarray}\label{efinite2}
E(m,x)w\in W[x,x^{-1}]\ \ \ \mbox{ for }m\in \Z,\; w\in W.
\end{eqnarray}}
\ed

It can be readily seen that indeed $C_{\rm fin}$ contains the
natural module $\C^{\infty}$.

The following is the main reason for formulating category $C_{\rm fin}$:

\bp{pchange} Let $W$ be a  $\widetilde{\gl}_{\infty}$-module of
level $\ell\in \C$ in category $C_{\rm fin}$. Then $W$ becomes a
restricted $\widetilde{\gl}_{\infty}^{e}$-module of level $\ell$
with
$$B(m,x)=E(m,e^{x})=\sum_{n\in \Z}E_{m,m+n}e^{-nx}$$
for $m\in \Z$. Furthermore, if $W$ is an irreducible
$\widetilde{\gl}_{\infty}$-module, then $W$ viewed as a
$\widetilde{\gl}_{\infty}^{e}$-module is also irreducible. On the
other hand, let $W_{1},W_{2}$ be $\widetilde{\gl}_{\infty}$-modules
in category $C_{\rm fin}$ and let $\theta$ be a linear map from
$W_{1}$ to $W_{2}$. Then $\theta$ is a homomorphism of
$\widetilde{\gl}_{\infty}^{e}$-modules if and only if $\theta$ is a
homomorphism of $\widetilde{\gl}_{\infty}$-modules. \ep

\begin{proof} First of all, notice that for any $p(x)\in \C[x,x^{-1}]$,
$p(e^{x})$ exists in $\C[[x]]$. {}From assumption we have
$E(m,x)w\in W[x,x^{-1}]$ for $m\in \Z,\; w\in W$. Then for each
$m\in \Z$, $E(m,e^{x})$ is a well defined element of $(\End
W)[[x]]$. For $m\in \Z$, set
$$\bar{B}(m,x)=E(m,e^{x})=\sum_{n\in \Z}E_{m,m+n}e^{-nx}.$$
Writing $\bar{B}(m,x)=\sum_{n\in \Z}\bar{B}(m,n)x^{-n-1}$, we have
$\bar{B}(m,r)=0$ for $r\ge 0$ and
\begin{eqnarray}
\bar{B}(m,-k-1)=\sum_{n\in \Z}\frac{1}{k!}(-n)^{k}E_{m,m+n}
\end{eqnarray}
for $k\ge 0$. Furthermore, we have
$$[\bar{B}(m,x_{1}),\bar{B}(n,x_{2})]
=e^{(m-n)(x_{1}-x_{2})}\left(\bar{B}(m,x_{2})-\bar{B}(n,x_{1})+f(m,n)\ell\right).$$
Thus $W$ becomes a $\widetilde{\gl}_{\infty}^{e}$-module of level
$\ell$ with $B(m,x)=\bar{B}(m,x)$ for $m\in \Z$, which is a
restricted module as $\bar{B}(m,x)$ involves only nonnegative powers
of $x$.

Regarding the assertion on irreducibility, let $w\in W,\; m\in \Z$.
For $k\ge 0$, we have
\begin{eqnarray}\label{equation-system}
(-1)^{k}k! \bar{B}(m,-k-1)w=\sum_{n\in \Z}n^{k}E_{m,m+n}w.
\end{eqnarray}
Note that the expression on the right hand side is a finite sum.
Considering $k\ge 1$, by solving the system of equations we can
express $E_{m,m+n}w$ for $n\ne 0$ in terms of $\bar{B}(m,-k-1)w$
with $k\ge 1$. We also have (with $k=0$)
$$\bar{B}(m,-1)w=\sum_{n\in \Z}E_{m,m+n}w.$$
Using this, we can express $E_{m,m}w$ in terms of $\bar{B}(m,-1)w$
and $E_{m,m+n}w$ for $n$ nonzero. Consequently, we can express
$E_{m,m+n}w$ for every $n\in \Z$ in terms of $\bar{B}(m,-k-1)w$ with
$k\ge 0$. It then follows immediately that $W$ viewed as a
$\widetilde{\gl}_{\infty}^{e}$-module is still irreducible.

The last assertion is also clear.
\end{proof}

In the following, we exhibit some irreducible
$\widetilde{\gl}_{\infty}$-modules of level $0$, or equivalently
$\gl_{\infty}$-modules, in category $C_{\rm fin}$. Note that
category $C_{\rm fin}$ is closed under direct sum and tensor
product. As $\C^{\infty}$ is in $C_{\rm fin}$, it follows that the
tensor algebra $T(\C^{\infty})$ is naturally a $\gl_{\infty}$-module
in $C_{\rm fin}$, on which $\gl_{\infty}$ acts by derivations.

\bex{examplesA} {\em Consider the symmetric algebra
$S(\C^{\infty})$, which is naturally a $\gl_{\infty}$-module in
$C_{\rm fin}$. Identify $S(\C^{\infty})$ with
polynomial algebra $\C[x_{n}\; |\ n\in \Z]$ on which $\gl_{\infty}$
acts by
$$E_{m,n}=x_{m}\frac{\partial}{\partial x_{n}}\ \ \ \mbox{ for
}m,n\in \Z.$$ Define $\deg x_{n}=1$ for $n\in \Z$ to make
$\C[x_{n}\; |\ n\in \Z]$ a $\Z$-graded algebra, and for $r\in \N$,
let $A_{r}$ denote the degree-$r$ homogeneous subspace, which is
linearly spanned by the monomials of total degree $r$. One can show
that $A_{r}$ ($r\in \Z$) are non-isomorphic irreducible submodules. } \eex

\bex{examples-ext} {\em On the other hand, the exterior algebra
$\Lambda(\C^{\infty})$ is also a $\gl_{\infty}$-module in $C_{\rm
fin}$, on which $\gl_{\infty}$ acts by derivations. Note that
$\Lambda(\C^{\infty})$ is an $\N$-graded algebra with $\C^{\infty}$
of degree $1$. Decompose $\Lambda(\C^{\infty})$ into homogeneous
subspaces
 $\Lambda(\C^{\infty})=\oplus_{r\in
\N}\Lambda^{r}(\C^{\infty})$. We have that
$\Lambda^{0}(\C^{\infty})=A_{0}=\C$ and
$\Lambda^{1}(\C^{\infty})=A_{1}=\C^{\infty}$. It is straightforward to show that
$\Lambda_{r}(\C^{\infty})$ with $r\in \Z$ are non-isomorphic
irreducible submodules.  Furthermore, one can show that $A_{m}$ (from Example \ref{examplesA}) and
$\Lambda^{n}(\C^{\infty})$ for $m,n\ge 2$ are non-isomorphic
irreducible $\gl_{\infty}$-modules.} \eex

\bex{examples-intermediate-series} {\em We here consider a generalization
of Example \ref{examplesA}. Let $S$ be any finite subset of $\Z$ and
let $\alpha: S\rightarrow \C$ be a function such that
$\alpha_{j}\notin \Z$ for $j\in S$. Set
$$V(S,\alpha)=\left(\prod_{j\in S}x_{j}^{\alpha_{j}}\right)
\C\left[x_{j},x_{j}^{-1}\ |\ j\in S\right]\otimes \C\left[x_{n}\ |\
n\in \Z\backslash S\right],$$ which is a $\gl_{\infty}$-module with
$E_{m,n}=x_{m}\frac{\partial}{\partial x_{n}}$ for $m,n\in \Z$. It
can be readily seen that $V(S,\alpha)$ is in category $C_{\rm fin}$.
Define $\deg x_{n}^{\lambda}=\lambda$ for $n\in \Z,\; \lambda\in
\C$, to make $V(S,\alpha)$ a $\C$-graded space. Then each
homogeneous subspace of $V(S,\alpha)$ is a submodule. One can show
that each homogeneous subspace as a $\gl_{\infty}$-module is
irreducible.} \eex

\br{rnon-example} {\em Note that $C_{\rm fin}$ does not contain any
nontrivial highest weight modules nor nontrivial lowest weight modules. On the
other hand, $C_{\rm fin}$ does not contain $\sl_{\infty}$ viewed as
an irreducible $\gl_{\infty}$-submodule of the adjoint module
either. All the examples given above are among what were called
intermediate series modules.} \er

It seems that every $\widetilde{\gl}_{\infty}$-module in category
$C_{\rm fin}$ is of level zero. Here we prove a weak version.

\bl{llevel=02} Let $W$ be a $\widetilde{\gl}_{\infty}$-module
satisfying the condition that for any $w\in W$, there exists a
finite subset $S$ of $\Z$ such that $E_{p,q}w=0$ for all $p,q\in \Z$
with $p,q\notin S$. Then $W$ is of level $0$. \el

\begin{proof} We may assume $W\ne 0$. Let $w$ be any nonzero vector in $W$.
By assumption, there exists a finite subset $S$ of $\Z$ such that
$E_{p,q}w=0$ for all $p,q\in \Z$ with $p,q\notin S$. As $S$ is
finite, there exist a negative integer $m$ and a positive integer
$n$ such that $m,n\notin S$. Then
$$E_{m,n}w=E_{n,m}w=E_{m,m}w=E_{n,n}w=0,$$
{}from which we get
$$0=E_{m,n}(E_{n,m}w)-E_{n,m}(E_{m,n}w)=E_{m,m}w-E_{n,n}w+\psi(E_{m,n},E_{n,m}){\bf k}w
={\bf k} w.$$ Therefore ${\bf k}$ acts trivially on $W$.
\end{proof}

\br{rvacuum-module} {\em Let $W$ be any
$\widetilde{\gl}_{\infty}$-module of level $\ell$ in category
$C_{\rm fin}$. By Proposition \ref{pchange}, $W$ becomes a
$\widetilde{\gl}_{\infty}^{e}$-module of level $\ell$. From the
construction, one sees that for any $m,n\in \Z$ with $n\ge 0$,
$B(m,n)$ acts on $W$ trivially. It follows from the construction of
$V_{\widetilde{\gl}_{\infty}^{e}}(\ell,0)$ that for any $w\in W$,
there exists a unique $\widetilde{\gl}_{\infty}^{e}$-module
homomorphism {}from $V_{\widetilde{\gl}_{\infty}^{e}}(\ell,0)$ to
$W$, sending ${\bf 1}$ to $w$.} \er

\bex{cinfinity} {\em From Example \ref{examplesA}, we have
non-isomorphic irreducible $\gl_{\infty}$-modules (or equivalently,
$\widetilde{\gl}_{\infty}$-modules of level $0$) $A_{r}$ ($r\in \N$)
in category $C_{\rm fin}$. By Proposition \ref{pchange}, we obtain
non-isomorphic irreducible $\widetilde{\gl}_{\infty}^{e}$-modules of
level $0$. The action of $\widetilde{\gl}_{\infty}^{e}$ on
$\C[x_{n}\ |\ n\in \Z]$ is determined by
$$B(m,x)x_{k}=x_{m}e^{(m-k)x}\ \ \ \mbox{ for }m,k\in \Z$$
with $B(m,n)$ ($m,n\in \Z$) acting as derivations. For $r\in \N$,
denote by $A_{r}^{e}$ the corresponding irreducible
$\widetilde{\gl}_{\infty}^{e}$-module. By Remark
\ref{rvacuum-module}, all irreducible modules $A_{r}^{e}$ for $r\in
\N$ are homomorphism images of
$V_{\widetilde{\gl}_{\infty}^{e}}(0,0)$. In view of this,
irreducible vacuum $\widetilde{\gl}_{\infty}^{e}$-modules of level
$0$ are not unique.} \eex

\br{rproblems} {\em There are two problems naturally arisen, one of which is
to classify irreducible $\widetilde{\gl}_{\infty}$-modules in category
$C_{\rm fin}$ and the other is to classify all irreducible vacuum
$\widetilde{\gl}_{\infty}^{e}$-modules of a fixed level $\ell$.
We hope to address these problems in a later publication.}
\er

We end this section with a generalized-Verma-module
construction. Recall that $\widetilde{\gl}_{\infty}^{e}[0]=E\otimes
\C[[t]]+ \C {\bf k}$ is an abelian subalgebra. Let $\ell\in \C$ and
let $\lambda: \Z\rightarrow \C$ be a function. Define a
$\widetilde{\gl}_{\infty}^{e}[0]$-module structure on $\C$ with
$E\otimes t\C[[t]]$ acting trivially, with $B(n,0)$ acting as
$\lambda_{n}$ for $n\in \Z$, and with ${\bf k}$ acting as $\ell$. We
denote this module by $\C_{\ell,\lambda}$. Then form an induced
module
\begin{eqnarray}
M(\ell,\lambda)=U(\widetilde{\gl}_{\infty}^{e})\otimes_{U(\widetilde{\gl}_{\infty}^{e}[0])}
\C_{\ell,\lambda}.
\end{eqnarray}
Noticing that
$\widetilde{\gl}_{\infty}^{e}=\widetilde{\gl}_{\infty}^{e}[0]\oplus
(E\otimes t^{-1}\C[t^{-1}])$ as a vector space, in view of the P-B-W
theorem we have
$$M(\ell,\lambda)=S(E\otimes t^{-1}\C[t^{-1}])$$
as a vector space.
 Notice that in case $\lambda=0$, we have
$M(\ell,0)=V_{\widetilde{\gl}_{\infty}^{e}}(\ell,0)$.

\br{rfinal} {\em Unlike the case for affine Lie algebras, it is not clear whether
$M(\ell,\lambda)$ has a unique maximal submodule. (Indeed, we have
seen that maximal submodules of
$M(\ell,0)\ \left(=V_{\widetilde{\gl}_{\infty}^{e}}(\ell,0)\right)$ (with $\lambda=0$) are not unique
as there are many non-isomorphic irreducible vacuum modules.)
It is also natural to ask under what condition $M(\ell,\lambda)$ is
irreducible. An immediate conjecture is that
$V_{\widetilde{\gl}_{\infty}^{e}}(\ell,0)$ is irreducible for any generic level $\ell$.}
\er


\begin{thebibliography}{AAGBP}
\bibitem[DJKM]{djkm}
E. Date, M. Jimbo, M. Kashiwara, T. Miwa, Operator approach to the
Kadomtsev-Petviashvili equation. Transformation groups for soliton
equations III, {\em J. Phys. Soc. Japan}, {\bf 50} (1981),
3806-3812.

\bibitem[DKM]{dkm}
E. Date, M. Kashiwara, T. Miwa, Vertex operators and $\tau$ functions. Transformation groups for soliton
equations II, {\em Proc. Japan Acad.} {\bf 57} Ser. A (1981),
387-392.

\bibitem[DL]{dl}
C. Dong and J. Lepowsky,  {\em Generalized Vertex Algebras and
Relative Vertex Operators}, Progress in Math., Vol. {\bf 112},
Birkh\"auser, Boston, 1993.

\bibitem[EK]{ek}
P. Etingof and D. Kazhdan, Quantization of Lie bialgebras, V, {\em
Selecta Mathematica, New Series}, {\bf 6} (2000) 105-130.

\bibitem[FZ]{fz}
I. Frenkel and Y.-C. Zhu, Vertex operator algebras associated to
representations of affine and Virasoro algebras, {\em Duke Math. J.}
{\bf 66} (1992) 123-168.

\bibitem[FKRW]{fkrw}
E. Frenkel, V. Kac, A. Radul, W. Wang, $\mathcal{W}_{1+\infty}$ and
$\mathcal{W}(\gl_{N})$ with central charge $N$,  {\em Commun. Math.
Phys.} {\bf 170} (1995) 337-357.

\bibitem[FLM]{flm}
I. B. Frenkel, J. Lepowsky and A. Meurman, {\em Vertex Operator
Algebras and the Monster,} Pure and Applied Math., Vol. 134,
Academic Press, Boston, 1988.

\bibitem[FZ]{fz}
I. B. Frenkel and Y. Zhu, Vertex operator algebras associated to
representations of affine and Virasoro algebras, {\em Duke Math. J.}
{\bf 66} (1992) 123-168.

\bibitem[JM]{jm}
M. Jimbo and T. Miwa, Solitons and infinite diemnsional Lie algebras, {\em Publ. RIMS, Kyoyo Univ.} {\bf 19} (1983) 943-1001.

\bibitem [K]{kac1}
V. Kac, {\it Infinite-dimensional Lie Algebras}, 3rd ed., Cambridge
Univ. Press, Cambridge, 1990.

\bibitem [KR]{kr}
V. Kac and A. Radul, Representation theory of the vertex algebra
$\mathcal{W}_{1+\infty}$, {\em Transf. Groups} {\bf 1} (1996) 41-70.

\bibitem[Li1]{li-local}
H.-S. Li, Local systems of vertex operators, vertex superalgebras
and modules,  {\em J. Pure Appl. Algebra} {\bf 109} (1996) 143-195.

\bibitem[Li2]{li-twisted}
H.-S. Li, Local systems of twisted vertex operators, vertex
superalgebras and twisted modules, in: {\em Moonshine, the Monster
and Related Topics, Proc. Joint Summer Research Conference, Mount
Holyoke,} 1994, ed. by C. Dong and G. Mason, Contemporary Math. {\bf
193}, Amer. Math. Soc., Providence, 1996, 203-236.

\bibitem[Li3]{li-qva1}
H.-S. Li, Nonlocal vertex algebras generated by formal vertex
operators, {\em Selecta Math. (New Series)} {\bf 11} (2005) 349-397.

\bibitem[Li4]{li-qva2}
H.-S. Li, Constructing quantum vertex algebras, {\em International
Journal of Mathematics} {\bf 17} (2006) 441-476.

\bibitem[Li5]{li-phi}
H.-S. Li, $\phi$-coordinated quasi modules for quantum vertex
algebras, {\em Commun. Math. Phys.} {\bf 308} (2011) 703-741.

\bibitem[Lia]{lian}
B. Lian, On the classification of simple vertex operator algebras,
{\em Commun. Math. Phys.} {\bf 164} (1994) 307-357.

\bibitem[Xu]{xu}
X.-P. Xu, Representations of centrally extended Lie algebras over
differential operators and vertex algebras, arXiv: math/0411146
[math.QA].

\end{thebibliography}
\end{document}